# CAN WE PROVE STABILITY BY USING A POSITIVE DEFINITE FUNCTION WITH NON SIGN-DEFINITE DERIVATIVE?

**Iasson Karafyllis**

**Department of Environmental Engineering, Technical University of Crete, 73100, Chania, Greece**

**Abstract**
Novel criteria for global asymptotic stability are presented. The results are obtained by a combination of the "discretization approach" and the ideas contained in the proof of the original Matrosov's result. The results can be used for the proof of global asymptotic stability by using continuously differentiable, positive definite functions which do not have a negative semi-definite derivative. Illustrating examples are provided.

## 1. Introduction

Lyapunov's direct method has been proved to be irreplaceable for the stability analysis of nonlinear systems. However, the main difficulty in the application of Lyapunov's direct method is to find a Lyapunov function for a given dynamical system. Most positive definite functions will not have a negative definite derivative for a given dynamical system and therefore cannot be used for stability analysis by using Lyapunov's direct method.

There are two ways to relax the requirement of a negative definite derivative:

1) By using the Krasovskii-La Salle principle (see [9,12,13,16]) or by using Matrosov's theorem (see [22]). The original result by Matrosov has been generalized recently in various directions (see [11,12,14,15,16,24]). However, in order to able to apply all available results it is necessary to have a positive definite (Lyapunov) function with negative semi-definite derivative or to assume uniform Lyapunov and Lagrange stability (which can be shown by a positive definite function with negative semi-definite derivative). It should be noted that the main idea in the proof of the original Matrosov's result is the division of the state space into two regions: in the first region (the "bad region") the non-positive derivative of the Lyapunov function can be arbitrarily small in absolute value while in the second region (the "good region") the derivative of the Lyapunov function has a negative upper bound. The proof is accomplished by showing that the solution cannot stay in the "bad region" forever and by estimating the time that the solution spends in the "good region". Recently, in [8] a different approach was proposed for a Lyapunov function which can have positive derivative in certain regions of the state space: by using the derivative of auxiliary functions the methodology guarantees that the solution enters the "good region" after a finite time and remains bounded. The idea of switching between system modes with negative and positive derivative of a Lyapunov function has been also used recently in the stability analysis of hybrid systems in [18,19]. However, in hybrid systems the time period for which the derivative of the Lyapunov function is positive is determined by the switching signal and it is not necessary to estimate it.

2) By using the "discretization approach" (see [1], the Appendix in [2] and recent generalizations in [17,20,21] as well as the proof of the main result in [6]), which does not require a negative definite derivative. Instead, the discretization approach requires that the difference of the Lyapunov function $V(x(T)) - V(x(0))$ is negative definite, where $x(t)$ denotes the solution of the dynamical system and $T > 0$ is a fixed time. Therefore, in this approach the Lyapunov function can even have a positive derivative in certain regions of the state space. However, the main difficulty in the application of this approach is the estimation of the difference of the Lyapunov function $V(x(T)) - V(x(0))$. The application of this approach to feedback stabilization problems gave very important results in [2] (see also recent extensions in [7]).

The purpose of the present work is to combine the above approaches and to provide useful global stability criteria that use a positive definite function with a non sign-definite derivative. The results are developed for the autonomous uncertain case

$$\dot{x} = f(d, x)$$
$$x \in \Re^n, d \in D \tag{1.1}$$

where $x(t)$ is the state and $d(t) \in D \subset \Re^l$ is a time-varying disturbance. However, the obtained results can be extended to the local case or the time-varying case. The key idea is the idea used in the proof of the original Matrosov's result described above concerning the division of the state space into two regions: the "good region", where the derivative of the Lyapunov function has a negative upper bound and the "bad region" where the derivative of the Lyapunov function can be positive. The first step is to show that the solution of (1.1) cannot stay in the "bad region" forever. Additional technical difficulties arise since we have to guarantee that the solution remains bounded while it stays in the "bad region". The second step is to estimate the difference of the Lyapunov function $V(x(T)) - V(x(0))$, where $T > 0$ is chosen appropriately so that the solution is in the "good region". Finally, by extending the discetization approach, we can guarantee robust global asymptotic stability or robust global exponential stability (Theorem 3.1, Corollary 3.5, Theorem 3.7 and Corollaries 3.8, 3.9, 3.10).

The structure of the paper is as follows. Section 2 provides the definitions of the notions used in the paper and some preliminary results that generalize the discretization approach. The results of Section 2 are interesting, since are necessary and sufficient conditions for robust global asymptotic stability. In Section 3 of the paper the main results are stated and proved. Illustrative examples of the proposed approach are provided in Section 4: the examples show how we can use very simple positive definite functions (e.g. $V(x) = |x|^2$), which do not have a sign-definite derivative. Finally, some concluding remarks are provided in Section 5. The Appendix contains the proofs of some technical steps needed in the proof of Theorem 3.1.

**Notations** Throughout this paper we adopt the following notations:
* For a vector $x \in \Re^n$ we denote by $|x|$ its usual Euclidean norm and by $x'$ its transpose.
* We say that an increasing continuous function $\gamma : \Re^+ \to \Re^+$ is of class $K$ if $\gamma(0) = 0$. We say that an increasing continuous function $\gamma : \Re^+ \to \Re^+$ is of class $K_\infty$ if $\gamma(0) = 0$ and $\lim_{s \to +\infty} \gamma(s) = +\infty$. By $KL$ we denote the set of all continuous functions $\sigma = \sigma(s,t) : \Re^+ \times \Re^+ \to \Re^+$ with the properties: (i) for each $t \geq 0$ the mapping $\sigma(\cdot, t)$ is of class $K$ ; (ii) for each $s \geq 0$, the mapping $\sigma(s, \cdot)$ is non-increasing with $\lim_{t \to +\infty} \sigma(s,t) = 0$.
* Let $D \subseteq \Re^l$ be a non-empty set. By $M_D$ we denote the class of all Lebesgue measurable and locally essentially bounded mappings $d : \Re^+ \to D$.
* By $C^j(A)$ ($C^j(A;\Omega)$), where $j \geq 0$ is a non-negative integer, $A \subseteq \Re^n$, we denote the class of functions (taking values in $\Omega \subseteq \Re^m$) that have continuous derivatives of order $j$ on $A$.
* For every scalar continuously differentiable function $V : \Re^n \to \Re$, $\nabla V(x)$ denotes the gradient of $V$ at $x \in \Re^n$, i.e., $\nabla V(x) = \left( \frac{\partial V}{\partial x_1}(x), \ldots, \frac{\partial V}{\partial x_n}(x) \right)$. We say that a function $V : \Re^n \to \Re^+$ is positive definite if $V(x) > 0$ for all $x \neq 0$ and $V(0) = 0$. We say that a continuous function $V : \Re^n \to \Re^+$ is radially unbounded if the following property holds: "for every $M > 0$ the set $\{ x \in \Re^n : V(x) \leq M \}$ is compact".

## 2. Preliminary Results

Throughout this paper we assume that system (1.1) satisfies the following hypotheses:

**(H1)** $D \subset \Re^l$ is compact.

**(H2)** The mapping $D \times \Re^n \ni (d,x) \to f(d,x) \in \Re^n$ is continuous with $f(d,0) = 0$ for all $d \in D$.

**(H3)** There exists a symmetric positive definite matrix $P \in \Re^{n \times n}$ such that for every compact set $S \subset \Re^n$ it holds that $\sup\left\{ \dfrac{(x-y)'P(f(d,x)-f(d,y))}{|x-y|^2} : d \in D, x, y \in S, x \neq y \right\} < +\infty$.

Hypothesis (H2) is a standard continuity hypothesis and hypothesis (H3) is often used in the literature instead of the usual local Lipschitz hypothesis for various purposes and is a generalization of the so-called "one-sided Lipschitz condition" (see, for example [23], page 416 and [3], page 106). Notice that the "one-sided Lipschitz condition" is weaker than the hypothesis of local Lipschitz continuity of the vector field $f(d,x)$ with respect to $x \in \Re^n$. It is clear that hypothesis (H3) guarantees that for every $(x_0, d) \in \Re^n \times M_D$, there exists a unique solution $x(t)$ of (1.1) with initial condition $x(0) = x_0$ corresponding to input $d \in M_D$. We denote by $x(t, x_0; d)$ the unique solution of (1.1) with initial condition $x(0) = x_0 \in \Re^n$ corresponding to input $d \in M_D$. Occasionally, we will use the following hypothesis for system (1.1):

**(H4)** For every compact set $S \subset \Re^n$ it holds that $\sup\left\{ \dfrac{|f(d,x)-f(d,y)|}{|x-y|} : d \in D, x, y \in S, x \neq y \right\} < +\infty$.

instead of hypothesis (H3). Hypothesis (H4) is more demanding than hypothesis (H3) in the sense that the implication (H4) $\Rightarrow$ (H3) holds.

We next continue by recalling the notion of Uniform (Robust) Global Asymptotic Stability.

**Definition 2.1:** *We say that $0 \in \Re^n$ is uniformly robustly globally asymptotically stable (URGAS) for system (1.1) under hypotheses (H1-3) if the following properties hold:*

- *for every $s > 0$, it holds that*

$$\sup\{|x(t, x_0; d)| ; t \geq 0, |x_0| \leq s, d \in M_D\} < +\infty$$
**(Uniform Robust Lagrange Stability)**

- *for every $\varepsilon > 0$ there exists a $\delta := \delta(\varepsilon) > 0$ such that:*
$$\sup\{|x(t, x_0; d)| ; t \geq 0, |x_0| \leq \delta, d \in M_D\} \leq \varepsilon$$
**(Uniform Robust Lyapunov Stability)**

- *for every $\varepsilon > 0$ and $s \geq 0$, there exists a $\tau := \tau(\varepsilon, s) \geq 0$, such that:*

$$\sup\{|x(t, x_0; d)| ; t \geq \tau, |x_0| \leq s, d \in M_D\} \leq \varepsilon$$
**(Uniform Attractivity for bounded sets of initial states)**

*For disturbance-free systems we say that $0 \in \Re^n$ is uniformly globally asymptotically stable (UGAS) for system (1.1).*

It should be noted that the notion of uniform robust global asymptotic stability coincides with the notion of uniform robust global asymptotic stability presented in [10]. We next provide the notion of global exponential stability (see also [9]).

**Definition 2.2:** *We say that $0 \in \Re^n$ is uniformly robustly globally exponentially stable (URGES) for (1.1) under hypotheses (H1-3), if there exist constants $M \geq 1$, $\sigma > 0$ such that the following inequality holds for all $t \geq 0$, $(x_0, d) \in \Re^n \times M_D$:*

$$|x(t, x_0; d)| \leq M \exp(-\sigma t)|x_0|$$

The following result is a generalization of the discretization approach for the autonomous case (1.1).

**Proposition 2.3:** *Consider system (1.1) under hypotheses (H1-3) and suppose that there exist a positive definite and radially unbounded $V \in C^0(\Re^n; \Re^+)$, a positive definite function $q \in C^0(\Re^+; \Re^+)$, a function $a \in K_\infty$ and a locally bounded function $T: \Re^n \setminus \{0\} \to (0, +\infty)$ such that for each $x_0 \in \Re^n \setminus \{0\}$, $d \in M_D$ the solution of (1.1) $x(t, x_0; d)$ with initial condition $x(0, x_0; d) = x_0$ corresponding to $d \in M_D$ exists on $[0, T(x_0)]$ and satisfies the following inequalities:*

$$V(x(t, x_0; d)) \le a(V(x_0)), \; \forall t \in [0, T(x_0)] \tag{2.1}$$

$$\min_{t \in [0, T(x_0)]} V(x(t, x_0; d)) \le V(x_0) - q(V(x_0)) \tag{2.2}$$

*Then $0 \in \Re^n$ is URGAS for (1.1). Moreover, if $T(x) \equiv r > 0$, $a(s) := Ms$, $q(s) := qs$, where $M, r > 0$, $q \in (0,1)$ and there exist constants $0 < K_1 < K_2$ with $K_1|x|^2 \le V(x) \le K_2|x|^2$ for all $x \in \Re^n$ then $0 \in \Re^n$ is Robustly Globally Exponentially Stable for (1.1).*

**Proof:** Let $\sigma \in KL$ be the function with the following property.

**Property (P):** if $\{V_i \ge 0\}_{i=0}^\infty$ is a sequence with $V_{i+1} \le V_i - q(V_i)$ then $V_i \le \sigma(V_0, i)$ for all $i \ge 0$.

The existence of $\sigma \in KL$ which satisfies property (P) is guaranteed by Lemma 4.3 in [4].

Define $T(0) = 1$. Let $x_0 \in \Re^n$, $d \in M_D$ arbitrary and define the following sequences:

$$x_{i+1} = x(t_i, x_i, P_{\tau_i} d), \; T_i = T(x_i), \; \tau_{i+1} = \tau_i + t_i, \; V_i = V(x_i), \; i \ge 0 \tag{2.3a}$$

with $\tau_0 = 0$, where $P_\tau d \in M_D$ is defined by $(P_\tau d)(t) = d(t + \tau)$ for $t \ge 0$ and $t_i \in [0, T_i]$ satisfies

$$V(x(t_i, x_i; P_{\tau_i} d)) = \min_{t \in [0, T_i]} V(x(t, x_i; P_{\tau_i} d)) \tag{2.3b}$$

for the case $x_i \ne 0$ and $t_i = T_i = 1$ for the case $x_i = 0$. Notice that by virtue of the semigroup property we obtain that $x_i = x(\tau_i, x_0; d)$.

Inequality (2.2) and definitions (2.3) imply that
$$V_{i+1} \le V_i - q(V_i) \tag{2.4}$$

for the case $x_i \ne 0$. For the case $x_i = 0$ by uniqueness of solution of (1.1) we have $x_{i+1} = 0$ and consequently inequality (2.4) holds as well in this case. Therefore, property (P) guarantees that

$$V_i \le \sigma(V_0, i) \text{ for all } i \ge 0 \tag{2.5}$$

where $\sigma \in KL$ is the function involved in property (P).

Inequality (2.1), definitions (2.3a) and the semigroup property guarantee that $V(x(t, x_0; d)) = V(x(t - \tau_i, x_i; P_{\tau_i} d)) \le a(V_i)$ for all $t \in [\tau_i, \tau_{i+1}]$ for the case $x_i \ne 0$. By uniqueness of solution of (1.1), it follows that $V(x(t, x_0; d)) = V(x(t - \tau_i, x_i; P_{\tau_i} d)) \le a(V_i)$ for all $t \in [\tau_i, \tau_{i+1}]$ for the case $x_i = 0$ as well. Since $\{V_i \ge 0\}_{i=0}^\infty$ is non-increasing (a consequence of (2.4)), we obtain

$$V(x(t, x_0; d)) \le a(V(x_0)) \text{ for all } t \in [0, \sup \tau_i) \tag{2.6}$$

Next we show that

$$V(x(t, x_0; d)) \le a(V(x_0)) \text{ for all } t \ge 0 \tag{2.7}$$

It should be noticed that Robust Lyapunov and Lagrange stability follows directly from inequality (2.7).

For the proof of inequality (2.7) we distinguish two cases:

<u>Case 1</u>: $\sup \tau_i < +\infty$. By virtue of inequality (2.5) we obtain that $\lim V_i = 0$ and consequently $\lim_{t \to (\sup \tau_i)^-} V(x(t, x_0; d)) = 0$. This implies that $\lim_{t \to (\sup \tau_i)^-} x(t, x_0; d) = 0$, which implies $x(t, x_0; d) = 0$ for all $t \geq \sup \tau_i$. Therefore inequality (2.7) is a consequence of (2.6) and the fact that $V(x(t, x_0; d)) = 0$ for all $t \geq \sup \tau_i$.

<u>Case 2</u>: $\sup \tau_i = +\infty$. In this case inequality (2.7) is a direct consequence of inequality (2.6).

We next show Robust Attractivity. Exploiting the fact that $V \in C^0(\Re^n; \Re^+)$ is a continuous, positive definite and radially unbounded function, it suffices to show that for every $\varepsilon > 0$, $R \geq 0$ there exists $\hat{T}(\varepsilon, R) \geq 0$ such that $V(x(t, x_0; d)) \leq \varepsilon$ for all $t \geq \hat{T}(\varepsilon, R)$, $x_0 \in \Re^n$ with $|x_0| \leq R$ and $d \in M_D$.

Let $\varepsilon > 0$, $R \geq 0$, $x_0 \in \Re^n$ with $|x_0| \leq R$ and $d \in M_D$ be arbitrary. By virtue of (2.7) and the semigroup property follows that if $V_i \leq a^{-1}(\varepsilon)$ for some $i \geq 0$ then we have $V(x(t, x_0; d)) \leq \varepsilon$ for all $t \geq \tau_i$. Define $J := \min\{i \geq 0 : V_i \leq a^{-1}(\varepsilon)\}$ and $B(R) := \max\{V(x) : |x| \leq R\}$. Let $N_\varepsilon(R) \in Z^+$ such that $\sigma(B(R), N_\varepsilon(R)) \leq a^{-1}(\varepsilon)$ and notice that inequality (2.5) and the fact that $V(x_0) \leq B(R)$ implies $J \leq N_\varepsilon(R)$.

Next suppose that $J \geq 1$. Since $V_i \geq a^{-1}(\varepsilon)$ for all $i \leq J-1$ (a consequence of definition $J := \min\{i \geq 0 : V_i \leq a^{-1}(\varepsilon)\}$), we get from (2.3a) and the facts that $\{V_i \geq 0\}_{i=0}^\infty$ is non-increasing (a consequence of (2.4)) and $V(x_0) \leq B(R)$:

$$\tau_{i+1} = \tau_i + t_i \leq \tau_i + T_i \leq \tau_i + \widetilde{T}(R) \text{ for all } i \leq J-1$$

where $\widetilde{T}_\varepsilon(R) := \sup\{T(x) : a^{-1}(\varepsilon) \leq V(x) \leq a^{-1}(\varepsilon) + B(R)\}$. Therefore $\tau_{i+1} \leq i\widetilde{T}_\varepsilon(R)$, for all $i \leq J-1$ and therefore inequality $J \leq N_\varepsilon(R)$ implies $\tau_J \leq N_\varepsilon(R)\widetilde{T}_\varepsilon(R)$. It follows that $V(x(t, x_0; d)) \leq \varepsilon$ for all $t \geq N_\varepsilon(R)\widetilde{T}_\varepsilon(R)$.

The above conclusion holds as well in the case $J = 0$, namely we have $V(x(t, x_0; d)) \leq \varepsilon$ for all $t \geq N_\varepsilon(R)\widetilde{T}_\varepsilon(R)$.

Thus for every $\varepsilon > 0$, $R \geq 0$ there exists $\hat{T}(\varepsilon, R) = N_\varepsilon(R)\widetilde{T}_\varepsilon(R) \geq 0$ such that $V(x(t, x_0; d)) \leq \varepsilon$ for all $t \geq \hat{T}(\varepsilon, R)$, $x_0 \in \Re^n$ with $|x_0| \leq R$ and $d \in M_D$.

Finally, for the case $T(x) \equiv r > 0$, $a(s) := Ms$, $q(s) := qs$, $K_1|x|^2 \leq V(x) \leq K_2|x|^2$ for all $x \in \Re^n$, where $M, r > 0$, $q \in (0,1)$ and $0 < K_1 < K_2$, we notice that inequality (2.4) implies $V_i \leq (1-q)^i V_0$ for all $i \geq 0$. Therefore, using the inequality $V(x(t, x_0; d)) \leq a(V_i)$ for all $t \in [\tau_i, \tau_{i+1}]$, in conjunction with definition $a(s) := Ms$ gives $V(x(t, x_0; d)) \leq (1-q)^i M V_0$ for all $t \in [\tau_i, \tau_{i+1}]$. The previous inequality combined with the inequalities $K_1|x|^2 \leq V(x) \leq K_2|x|^2$ for all $x \in \Re^n$ gives $|x(t, x_0; d)| \leq \exp(-\sigma i)\sqrt{\frac{K_2 M}{K_1}}|x_0|$, for all $t \in [\tau_i, \tau_{i+1}]$, where $\sigma > 0$ is defined by the equation $\exp(-2\sigma) = 1-q$. Using the fact that $t \leq \tau_{i+1} \leq (i+1)r$, we obtain the inequality $\exp(-\sigma i) \leq \exp(\sigma)\exp\left(-\frac{\sigma}{r}t\right)$, for all $i \geq 0$ and $t \in [\tau_i, \tau_{i+1}]$. Consequently, by distinguishing again the cases $\sup \tau_i < +\infty$ and $\sup \tau_i = +\infty$, we have $|x(t, x_0; d)| \leq \exp\left(-\frac{\sigma}{r}t\right)\exp(\sigma)\sqrt{\frac{K_2 M}{K_1}}|x_0|$, for all $t \geq 0$, which implies that $0 \in \Re^n$ is Robustly Globally Exponentially Stable for (1.1).

The proof is complete. ◁

**Remark 2.4:** The reader should notice that the converse of Proposition 2.3 holds, i.e., if $0 \in \Re^n$ is URGAS for (1.1) then for every positive definite and radially unbounded function $V \in C^0(\Re^n; \Re^+)$, there exist a function $a \in K_\infty$ and a locally bounded function $T : \Re^n \setminus \{0\} \to (0, +\infty)$ such that for each $x_0 \in \Re^n \setminus \{0\}$, $d \in M_D$ the solution of (1.1) $x(t, x_0; d)$ with initial condition $x(0, x_0; d) = x_0$ corresponding to $d \in M_D$ exists on $[0, T(x_0)]$ and satisfies the inequalities (2.1) and (2.2). Indeed, if $0 \in \Re^n$ is URGAS for (1.1) then there exists $\sigma \in KL$ such that for each $x_0 \in \Re^n$, $d \in M_D$ the solution of (1.1) $x(t, x_0; d)$ with initial condition $x(0, x_0; d) = x_0$ corresponding to $d \in M_D$ satisfies:

$$|x(t, x_0; d)| \leq \sigma(|x_0|, t), \quad \forall t \geq 0 \tag{2.8}$$

Without loss of generality, we may assume that for each $s > 0$ the mapping $t \to \sigma(s, t)$ is strictly decreasing (if not replace $\sigma(s, t)$ by $\sigma(s, t) + s \exp(-t)$). Since $V \in C^0(\Re^n; \Re^+)$ is positive definite and radially unbounded, there exist functions $a_1, a_2 \in K_\infty$ such that

$$a_1(|x|) \leq V(x) \leq a_2(|x|), \quad \forall x \in \Re^n \tag{2.9}$$

Combining (2.8) and (2.9) we obtain:

$$V(x(t, x_0; d)) \leq a_2\left(\sigma\left(a_1^{-1}(V(x_0)), t\right)\right), \quad \forall t \geq 0 \tag{2.10}$$

Let $q \in (0,1)$ and let $t(s) > 0$ be the solution of the equation $a_2\left(\sigma\left(a_1^{-1}(s), t(s)\right)\right) = (1-q)s$, for each $s > 0$. It can be shown by contradiction that the mapping $(0, +\infty) \ni s \to t(s)$ is bounded on every compact set $S \subset (0, +\infty)$. Therefore, by virtue of (2.10), we conclude that inequalities (2.1) and (2.2) hold with $a(s) := a_2\left(\sigma\left(a_1^{-1}(s), 0\right)\right)$, $q(s) := q s$ and $T(x) := t(V(x)) + 1$. ◁

The following proposition is less demanding in terms of the inequalities that guarantee URGAS. However, in contrast to Proposition 2.3, we have to assume that system (1.1) is forward complete and that hypothesis (H4) holds. We say that system (1.1) is forward complete if for every $x_0 \in \Re^n$, $d \in M_D$ the solution of (1.1) $x(t, x_0; d)$ with initial condition $x(0, x_0; d) = x_0$ corresponding to $d \in M_D$ is defined for all $t \geq 0$.

**Proposition 2.5:** *Consider system (1.1) under hypotheses (H1), (H2), (H4) and assume that system (1.1) is forward complete. Furthermore, suppose that there exist a positive definite and radially unbounded $V \in C^0(\Re^n; \Re^+)$, a positive definite function $q \in C^0(\Re^+; \Re^+)$ and a locally bounded function $T : \Re^n \to (0, +\infty)$ such that for each $x_0 \in \Re^n \setminus \{0\}$, $d \in M_D$ the solution of (1.1) $x(t, x_0; d)$ with initial condition $x(0, x_0; d) = x_0$ corresponding to $d \in M_D$ exists on $[0, T(x_0)]$ and satisfies inequality (2.2). Then $0 \in \Re^n$ is URGAS for (1.1).*

The reader should notice that an additional difference between Proposition 2.5 and Proposition 2.3 is the fact that Proposition 2.3 demands the function $T : \Re^n \setminus \{0\} \to (0, +\infty)$ to be locally bounded while Proposition 2.5 demands the function $T : \Re^n \to (0, +\infty)$ to be locally bounded. Since the value $T(0)$ plays no role, it is clear that the extra assumption required for Proposition 2.5 can be replaced by the condition $\limsup_{x \to 0} T(x) < +\infty$.

**Proof of Proposition 2.5:** The key idea of the proof is to show that forward completeness + hypothesis (H4) + $\limsup_{x \to 0} T(x) < +\infty$ imply the existence of a function $a \in K_\infty$ such that inequality (2.1) holds as well. Then Proposition 2.3 guarantees that $0 \in \Re^n$ is URGAS for (1.1).

Indeed, since (1.1) is forward complete and since hypothesis (H4) holds, Proposition 5.1 in [5] guarantees that system (1.1) is Robustly Forward Complete (RFC, see [5]). Lemma 2.3 in [5] guarantees the existence of functions $\mu \in K^+$ and $\zeta \in K_\infty$ such that the following inequality holds for all $x_0 \in \Re^n$, $d \in M_D$ and $t \geq 0$:

$$|x(t, x_0; d)| \leq \mu(t)\zeta(|x_0|) \tag{2.11}$$

Without loss of generality, we may assume that $\mu \in K^+$ is non-decreasing. Since $V \in C^0(\Re^n; \Re^+)$ is positive definite and radially unbounded, there exist functions $a_1, a_2 \in K_\infty$ such that inequality (2.9) holds. Combining (2.9) and (2.11) we obtain that for each $x_0 \in \Re^n \setminus \{0\}$, $d \in M_D$ the solution $x(t, x_0; d)$ of (1.1) satisfies the following inequality:

$$V(x(t, x_0; d)) \leq a_2\big(\mu(T(x_0))\zeta(|x_0|)\big), \quad \forall t \in [0, T(x_0)] \tag{2.12}$$

Define $p(s) := a_2\left(\mu\left(\sup_{|x| \leq s} T(x)\right)\zeta(s)\right)$ for all $s \geq 0$. Since $T: \Re^n \to (0, +\infty)$ is locally bounded, it follows that $p(s)$ is well-defined for all $s \geq 0$ and is a non-decreasing function. Moreover, it holds that $\lim_{s \to 0^+} p(s) = p(0) = 0$. Define $\tilde{a}(s) := s + \frac{1}{s}\int_s^{2s} p(\xi)d\xi$ for $s > 0$ and $\tilde{a}(0) := 0$. The function $\tilde{a}$ is of class $K_\infty$ and satisfies $\tilde{a}(s) \geq p(s)$ for all $s \geq 0$. Consequently, using (2.12) we obtain that for each $x_0 \in \Re^n \setminus \{0\}$, $d \in M_D$ the solution $x(t, x_0; d)$ of (1.1) satisfies the following inequality:

$$V(x(t, x_0; d)) \leq \tilde{a}(|x_0|), \quad \forall t \in [0, T(x_0)] \tag{2.13}$$

Using (2.9) and (2.13), it follows that inequality (2.1) holds with $a(s) := \tilde{a}(a_1^{-1}(s))$. The proof is complete. ◁

## 3. Main Results

Next, the main result of the present work is stated.

**Theorem 3.1:** *Consider system (1.1) under hypotheses (H1), (H2), (H4) and suppose that there exist a positive definite and radially unbounded function $V \in C^1(\Re^n; \Re^+)$, a family of functions $W_i \in C^1(\Re^n; \Re)$ with $W_i(0) = 0$, $b_i \in C^0(\Re^+; \Re^+)$ ($i = 0, ..., k$), $\rho, c_1, c_2, g, \gamma \in K_\infty$, with $\rho(s) > c_1(s) \geq c_2(s)$ for all $s > 0$, $\lambda \in K$ with $\lambda(s) < s$ for all $s > 0$, a locally bounded function $r: \Re^+ \to (0, +\infty)$ and a $C^1$ function $\mu: \Re^+ \to \Re$ with $\mu(0) = 0$ for which the function $\kappa(s) := c_1(s) + \mu(s)$ is non-decreasing, such that the following inequalities hold:*

$$\max_{d \in D} \nabla V(x) f(d, x) \leq -\rho(V(x)) + W_0(x), \text{ for all } x \in \Re^n \tag{3.1}$$

$$\max_{d \in D} \nabla W_i(x) f(d, x) \leq W_{i+1}(x), \text{ for all } i = 0, ..., k-1 \text{ and for all } x \in \Re^n \text{ with } W_0(x) \geq c_2(V(x)) \tag{3.2}$$

$$W_i(x) \leq b_i(V(x)), \text{ for all } i = 0, ..., k \text{ and for all } x \in \Re^n \text{ with } W_0(x) \geq c_2(V(x)) \tag{3.3}$$

$$\max_{d \in D} \nabla W_k(x) f(d, x) \leq -g(V(x)), \text{ for all } x \in \Re^n \text{ with } W_0(x) \geq c_2(V(x)) \tag{3.4}$$

$$\max_{d \in D} \nabla W_0(x) f(d, x) + \max_{d \in D} \mu'(V(x))\nabla V(x) f(d, x) \leq 0, \text{ for all } x \in \Re^n \text{ with } c_1(V(x)) \geq W_0(x) \geq c_2(V(x)) \tag{3.5}$$

$$c_1(\lambda(s)) + \mu(\lambda(s)) > c_2(\gamma(s)) + \mu(\gamma(s)), \text{ for all } s > 0 \tag{3.6}$$

$$c_2(\lambda(s)) + g(\lambda(s))\frac{r^{k+1}(s)}{(k+1)!} > \sum_{i=0}^{k} \frac{r^i(s)}{i!} b_i(s), \text{ for all } s > 0 \qquad (3.7)$$

$$\gamma(s) \geq \max\left\{ s, \rho^{-1}\left( \max_{\tau \in [0, r(s)]} \left[ \sum_{i=0}^{k} \frac{\tau^i}{i!} b_i(s) - g(\lambda(s)) \frac{\tau^{k+1}}{(k+1)!} \right] \right) \right\}, \text{ for all } s > 0 \qquad (3.8)$$

$$\limsup_{s \to 0^+} \int_{\lambda(s)}^{\gamma(s)} \frac{d\tau}{\rho(\tau) - c_1(\tau)} < +\infty \qquad (3.9)$$

Then $0 \in \Re^n$ is URGAS for (1.1).

**Remark 3.2:** A sufficient condition for the existence of a locally bounded function $r : \Re^+ \to (0, +\infty)$ that satisfies (3.7) is the set of inequalities $\limsup_{s \to 0^+} \frac{b_i(s)}{g(\lambda(s))} < +\infty$, for $i = 1, \ldots, k$ and $\limsup_{s \to 0^+} \frac{b_0(s) - c_2(\lambda(s))}{g(\lambda(s))} < +\infty$. More specifically, if the previous set of inequalities holds then the map defined by

$$r(s) := 1 + \max\left\{ \left( (k+1)! \frac{(k+1)(b_0(s) - c_2(\lambda(s)))}{g(\lambda(s))} \right)^{\frac{1}{k+1}}, \max_{i=1, \ldots, k} \left( \frac{(k+1)!}{i!} \frac{(k+1) b_i(s)}{g(\lambda(s))} \right)^{\frac{1}{k+1-i}} \right\} \text{ for } s > 0 \text{ and } r(0) := 1, \text{ is a}$$

locally bounded function $r : \Re^+ \to (0, +\infty)$ that satisfies (3.7).

**Remark 3.3:** A sufficient condition for (3.9) is the existence of a constant $K \in (0,1)$ such that:

$$\rho(s) \geq c_1(s) + K s \text{ and } \gamma(s) \leq K^{-1} \lambda(s), \text{ for } s > 0 \text{ sufficiently small}$$

However, the above condition is not necessary (e.g., if $\rho(s) \geq c_1(s) + K s^2$ and $\lambda(s) \geq \frac{\gamma(s)}{1 + \gamma(s)}$ for $s > 0$ sufficiently small then (3.9) holds).

The proof of Theorem 3.1 is heavily based on the following lemma. Its proof can be found in the Appendix.

**Lemma 3.4:** *Consider system (1.1) under hypotheses (H1), (H2), (H4) and suppose that there exist a positive definite and radially unbounded function $V \in C^1(\Re^n; \Re^+)$, a family of functions $W_i \in C^1(\Re^n; \Re)$ with $W_i(0) = 0$, $b_i \in C^0(\Re^+; \Re^+)$ ($i = 0, \ldots, k$), $\rho, c_1, c_2, g, \gamma \in K_\infty$, with $\rho(s) > c_1(s) \geq c_2(s)$ for all $s > 0$, $\lambda \in K$ with $\lambda(s) < s$ for all $s > 0$, a locally bounded function $r : \Re^+ \to (0, +\infty)$ and a $C^1$ function $\mu : \Re^+ \to \Re$ with $\mu(0) = 0$ for which the function $\kappa(s) := c_1(s) + \mu(s)$ is non-decreasing, such that inequalities hold (3.1)-(3.8) hold. Then system (1.1) is forward complete.*

We are now ready to provide the proof of Theorem 3.1.

**Proof of Theorem 3.1:** By virtue of Proposition 2.5 and Lemma 3.4, it suffices to show that inequalities (2.1), (2.2) hold for each $x_0 \in \Re^n \setminus \{0\}$, $d \in M_D$ for the function $V \in C^1(\Re^n; \Re^+)$ with

$$T(x) := p(V(x)) \qquad (3.10)$$

$$p(s) := r(s) + \int_{\lambda(s)}^{\gamma(s)} \frac{d\tau}{\rho(\tau) - c_1(\tau)}, \text{ for } s > 0 \text{ and } p(0) := 1 \qquad (3.11)$$

$$q(s) := s - \lambda(s) \tag{3.12}$$

The reader should notice that condition (3.9) and the fact that $r: \Re^+ \to (0, +\infty)$ is a locally bounded function guarantee that the function $p: \Re^+ \to (0, +\infty)$ as defined by (3.11) is locally bounded.

Let $x_0 \in \Re^n \setminus \{0\}$, $d \in M_D$ (arbitrary). We next show by contradiction that there exists $t \in [0, T(x_0)]$ such that $V(x(t, x_0; d)) \leq \lambda(V(x_0))$. Then definition (3.12) automatically guarantees that inequality (2.2) holds.

Assume next that $V(x(t, x_0; d)) > \lambda(V(x_0))$ for all $t \in [0, T(x_0)]$. We show that this cannot happen.

We first start by stating the following fact. Its proof can be found in the Appendix.

FACT I: There exists $t \in [0, r(V(x_0))]$ with $W_0(x(t, x_0; d)) < c_2(V(x(t, x_0; d)))$.

Define:

$$t_1 := \inf\{t \in [0, r(V(x_0))]: W_0(x(t, x_0; d)) < c_2(V(x(t, x_0; d)))\} \tag{3.13}$$

We next continue with the following fact. Its proof can be found in the Appendix.

FACT II: The following inequality holds:

$$V(x(t_1, x_0; d)) \leq \gamma(V(x_0)) \tag{3.14}$$

We next distinguish the following cases:

CASE 1: $W_0(x(t, x_0; d)) \leq c_1(V(x(t, x_0; d)))$ for all $t \in [t_1, T(x_0)]$.

Define $c(s) := \rho(s) - c_1(s)$, which is a positive definite, continuous function. In this case inequality (3.1) implies that $\dot{V}(t) \leq -c(V(t))$, for $t \in [t_1, T(x_0)]$, a.e., where $V(t) := V(x(t, x_0; d))$. Consequently, we obtain $\int_{V(t_1)}^{V(T(x_0))} \frac{ds}{c(s)} = \int_{t_1}^{T(x_0)} \frac{\dot{V}(t)}{c(V(t))} dt \leq -(T(x_0) - t_1)$. Combining, the previous inequality with the fact that $t_1 \in [0, r(V(x_0))]$ and definition (3.11) we get $\int_{V(T(x_0))}^{V(t_1)} \frac{ds}{c(s)} \geq \int_{\lambda(V(x_0))}^{\gamma(V(x_0))} \frac{ds}{c(s)}$. Since (3.14) holds, the previous inequality gives $V(T(x_0)) \leq \lambda(V(x_0))$, a contradiction.

CASE 2: There exists $t \in [t_1, T(x_0)]$ with $W_0(x(t, x_0; d)) > c_1(V(x(t, x_0; d)))$.

In this case, continuity of mapping $t \to \frac{W_0(x(t, x_0; d))}{V(x(t, x_0; d))}$ guarantees the existence of times $t_2 < t_3$ with $t_1 \leq t_2 < t_3 \leq T(x_0)$ and such that:

$$W_0(x(t, x_0; d)) \leq c_1(V(x(t, x_0; d))), \text{ for all } t \in [t_1, t_3] \tag{3.15}$$

$$W_0(x(t_2, x_0; d)) = c_2(V(x(t_2, x_0; d))), \quad W_0(x(t_3, x_0; d)) = c_1(V(x(t_3, x_0; d))) \tag{3.16}$$

$$W_0(x(t, x_0; d)) \geq c_2(V(x(t, x_0; d))), \text{ for all } t \in [t_2, t_3] \quad (3.17)$$

Inequality (3.15) in conjunction with inequality (3.1) guarantees that

$$V(t) \leq V(t_2) - \int_{t_2}^{t} (\rho(V(\tau)) - c_1(V(\tau))) d\tau, \text{ for all } t \in [t_2, t_3] \quad (3.18)$$

$$V(t) \leq V(t_1), \text{ for all } t \in [t_1, t_3] \quad (3.19)$$

Inequalities (3.15), (3.17) and (3.5) imply that

$$W_0(t_3) + \mu(V(t_3)) \leq W_0(t_2) + \mu(V(t_2)) \quad (3.20)$$

It follows from (3.16) and (3.20) that:

$$c_1(V(t_3)) + \mu(V(t_3)) \leq c_2(V(t_2)) + \mu(V(t_2)) \quad (3.21)$$

If $V(t_3) \leq \lambda(\gamma^{-1}(V(t_2)))$ then using (3.14), (3.19) we obtain $V(t_3) \leq \lambda(V(x_0))$, a contradiction.

Thus we are left with the case $V(t_3) > \lambda(\gamma^{-1}(V(t_2)))$. In this case inequality (3.21) and the fact that the function $\kappa(s) := c_1(s) + \mu(s)$ is non-decreasing, give

$$c_1(\lambda(\gamma^{-1}(V(t_2)))) + \mu(\lambda(\gamma^{-1}(V(t_2)))) \leq c_2(V(t_2)) + \mu(V(t_2))$$

The above inequality contradicts inequality (3.6) for $s = \gamma^{-1}(V(t_2))$.

The proof is complete. ◁

**Corollary 3.5:** *Consider system (1.1) under hypotheses (H1), (H2), (H4) and suppose that there exist a positive definite and radially unbounded function $V \in C^1(\Re^n; \Re^+)$, a locally Lipschitz function $\phi : \Re^n \to (0, +\infty)$, a family of functions $W_i \in C^1(\Re^n; \Re)$ with $W_i(0) = 0$, $b_i \in C^0(\Re^+; \Re^+)$ ($i = 0, ..., k$), $\rho, c_1, c_2, g, \gamma \in K_\infty$, with $\rho(s) > c_1(s) \geq c_2(s)$ for all $s > 0$, $\lambda \in K$ with $\lambda(s) < s$ for all $s > 0$, a locally bounded function $r : \Re^+ \to (0, +\infty)$ and a $C^1$ function $\mu : \Re^+ \to \Re$ with $\mu(0) = 0$ for which the function $\kappa(s) := c_1(s) + \mu(s)$ is non-decreasing, such that inequalities (3.3), (3.5), (3.6), (3.7), (3.8), (3.9) hold as well as the following inequalities:*

$$\max_{d \in D} \nabla V(x) f(d, x) \leq -\phi(x) \rho(V(x)) + \phi(x) W_0(x), \text{ for all } x \in \Re^n \quad (3.22)$$

$$\max_{d \in D} \nabla W_i(x) f(d, x) \leq \phi(x) W_{i+1}(x), \text{ for all } i = 0, ..., k-1 \text{ and for all } x \in \Re^n \text{ with } W_0(x) \geq c_2(V(x)) \quad (3.23)$$

$$\max_{d \in D} \nabla W_k(x) f(d, x) \leq -\phi(x) g(V(x)), \text{ for all } x \in \Re^n \text{ with } W_0(x) \geq c_2(V(x)) \quad (3.24)$$

*Then $0 \in \Re^n$ is URGAS for (1.1).*

**Proof:** Simply consider the dynamical system:

$$\dot{x} = \frac{1}{\phi(x)} f(d, x)$$
$$x \in \Re^n, d \in D \quad (3.25)$$

Since $\phi: \Re^n \to (0,+\infty)$ is locally Lipschitz, it follows that system (3.25) satisfies hypotheses (H1), (H2) and (H4). Moreover, all requirements of Theorem 3.1 are fulfilled and consequently $0 \in \Re^n$ is URGAS for (3.25). Classical Lyapunov theory implies that $0 \in \Re^n$ is URGAS for (1.1). The proof is complete. ◁

**Remark 3.6:** Here it should be noticed that Lyapunov's direct method is a special case of Corollary 3.5. Indeed, if there exists a positive definite continuous function $q: \Re^n \to \Re^+$ such that $\max_{d \in D} \nabla V(x) f(d,x) \leq -q(x)$ for all $x \in \Re^n$, then one can construct a locally Lipschitz function $\phi: \Re^n \to (0,+\infty)$ and a function $\rho \in K_\infty$ such that $\max_{d \in D} \nabla V(x) f(d,x) \leq -\phi(x)\rho(V(x))$ for all $x \in \Re^n$. Consequently, inequality (3.22) holds with $W_0(x) \equiv 0$. Therefore, all requirements of Corollary 3.5 are satisfied with $k = 0$. The reader should notice that since $k = 0$ inequalities (3.23) do not apply and since $W_0(x) \equiv 0$, the set of all $x \in \Re^n$ with $W_0(x) \geq c_2(V(x))$ is reduced to the singleton $\{0\}$ for every $c_2 \in K_\infty$. Hence, inequality (3.24) holds with arbitrary $g \in K_\infty$. Moreover, inequality (3.3) holds with $b_0(s) \equiv 0$, inequality (3.5) holds with $\mu(s) \equiv 0$ and inequalities (3.6), (3.7), (3.8), (3.9) hold with $r(s) \equiv 1$, $\gamma(s) := s$, $c_1(s) := \frac{3}{4}\rho(s)$, $c_2(s) := \frac{1}{2}\rho\left(\frac{s}{2}\right)$ and arbitrary $\lambda \in K$ with $\lambda(s) < s$ for all $s > 0$, which satisfies $\lambda(s) \geq \max\left\{\frac{s}{2}, s - \frac{1}{4}\rho\left(\frac{s}{2}\right)\right\}$ for $s > 0$ sufficiently small.

The following theorem provides stability criteria under minimal regularity requirements for system (1.1). Here we do not assume the local Lipschitz assumption (H4).

**Theorem 3.7:** *Consider system (1.1) under hypotheses (H1), (H2), (H3) and suppose that there exist a positive definite and radially unbounded function $V \in C^1(\Re^n; \Re^+)$, a family of functions $W_i \in C^1(\Re^n; \Re)$ with $W_i(0) = 0$, constants $b_i \geq 0$ ($i = 0,...,k$) with $b_0 \geq \rho$, $\rho, c_1, c_2, g, \gamma, r > 0$, with $\rho > c_1 \geq c_2 > 0$, $\lambda \in (0,1)$ and $\mu \geq -c_1$ such that the following inequalities hold:*

$$\max_{d \in D} \nabla V(x) f(d,x) \leq -\rho V(x) + W_0(x), \text{ for all } x \in \Re^n \tag{3.26}$$

$$\max_{d \in D} \nabla W_i(x) f(d,x) \leq W_{i+1}(x), \text{ for all } i = 0,...,k-1 \text{ and for all } x \in \Re^n \text{ with } W_0(x) \geq c_2 V(x) \tag{3.27}$$

$$W_i(x) \leq b_i V(x), \text{ for all } i = 0,...,k \text{ and for all } x \in \Re^n \text{ with } W_0(x) \geq c_2 V(x) \tag{3.28}$$

$$\max_{d \in D} \nabla W_k(x) f(d,x) \leq -g V(x), \text{ for all } x \in \Re^n \text{ with } W_0(x) \geq c_2 V(x) \tag{3.29}$$

$$\max_{d \in D} \nabla W_0(x) f(d,x) + \max_{d \in D} \mu \nabla V(x) f(d,x) \leq 0, \text{ for all } x \in \Re^n \text{ with } c_1 V(x) \geq W_0(x) \geq c_2 V(x) \tag{3.30}$$

$$(c_1 + \mu)\lambda > (c_2 + \mu)\gamma \tag{3.31}$$

$$c_2 \lambda + g\lambda \frac{r^{k+1}}{(k+1)!} > \sum_{i=0}^{k} \frac{r^i}{i!} b_i \tag{3.32}$$

$$\gamma \geq \min\left\{\exp((b_0 - \rho)r); \max_{\tau \in [0,r]} \frac{1}{\rho}\left[\sum_{i=0}^{k} \frac{\tau^i}{i!} b_i - g\lambda \frac{\tau^{k+1}}{(k+1)!}\right]\right\} \tag{3.33}$$

*Then $0 \in \Re^n$ is URGAS for (1.1). Moreover, if there exist constants $0 < K_1 < K_2$ with $K_1|x|^2 \leq V(x) \leq K_2|x|^2$ for all $x \in \Re^n$ then $0 \in \Re^n$ is Robustly Globally Exponentially Stable for (1.1).*

**Proof:** Let $x_0 \in \Re^n \setminus \{0\}$ and $d \in M_D$ (arbitrary). Inequalities (3.26), (3.28) (for $i = 0$) and the fact $b_0 \geq \rho$ imply that

$$V(x(t, x_0; d)) \leq \exp((b_0 - \rho)t)V(x_0), \text{ for all } t \geq 0 \tag{3.34}$$

Indeed, inequality (3.26) implies $\max_{d \in D} \nabla V(x) f(d, x) \leq 0$, when $W_0(x) \leq c_2 V(x)$. Moreover, inequalities (3.26), (3.28) (for $i = 0$) imply that $\max_{d \in D} \nabla V(x) f(d, x) \leq (b_0 - \rho)V(x)$, when $W_0(x) \geq c_2 V(x)$. Since $b_0 \geq \rho$, we conclude that $\max_{d \in D} \nabla V(x) f(d, x) \leq (b_0 - \rho)V(x)$ for all $x \in \Re^n$. Inequality (3.34) follows directly from the previous differential inequality.

The proof is exactly the same with the proof of Theorem 3.1 with $\rho(s) := \rho s$, $c_1(s) := c_1 s$, $c_2(s) := c_2 s$, $g(s) := g s$, $\gamma(s) := \gamma s$, $\lambda(s) := \lambda s$, $r(s) \equiv r$, $\mu(s) := \mu s$, with two major differences:

1) Instead of working with Proposition 2.5, we are working with Proposition 2.3. Indeed, inequality (3.34) guarantees that inequality (2.1) holds with $a(s) := s \exp((b_0 - \rho)T)$, where $T := r + \dfrac{\ln(\gamma) - \ln(\lambda)}{(\rho - c_1)}$.

2) Inequality (3.14) is obtained by a combined use of the proof of Fact II in the Appendix and inequality (3.34).

Details are left to the reader. ◁

**Corollary 3.8:** *Consider system (1.1) under hypotheses (H1), (H2), (H3) and suppose that there exist a positive definite and radially unbounded function $V \in C^1(\Re^n; \Re^+)$, a locally Lipschitz function $\phi: \Re^n \to (0, +\infty)$, a family of functions $W_i \in C^1(\Re^n; \Re)$ with $W_i(0) = 0$, constants $b_i \geq 0$ ($i = 0, ..., k$) with $b_0 \geq \rho$, $\rho, c_1, c_2, g, \gamma, r > 0$, with $\rho > c_1 \geq c_2 > 0$, $\lambda \in (0,1)$ and $\mu \geq -c_1$ such that inequalities (3.28), (3.30), (3.31), (3.32), (3.33) as well as the following inequalities hold:*

$$\max_{d \in D} \nabla V(x) f(d, x) \leq -\rho \phi(x) V(x) + \phi(x) W_0(x), \text{ for all } x \in \Re^n \tag{3.35}$$

$$\max_{d \in D} \nabla W_i(x) f(d, x) \leq \phi(x) W_{i+1}(x), \text{ for all } i = 0, ..., k-1 \text{ and for all } x \in \Re^n \text{ with } W_0(x) \geq c_2 V(x) \tag{3.36}$$

$$\max_{d \in D} \nabla W_k(x) f(d, x) \leq -g \phi(x) V(x), \text{ for all } x \in \Re^n \text{ with } W_0(x) \geq c_2 V(x) \tag{3.37}$$

*Then $0 \in \Re^n$ is URGAS for (1.1). Moreover, if there exist constants $0 < K_1 < K_2$ with $K_1 |x|^2 \leq V(x) \leq K_2 |x|^2$ and $\phi(x) \geq K_1$ for all $x \in \Re^n$ then $0 \in \Re^n$ is Robustly Globally Exponentially Stable for (1.1).*

**Proof:** Again the proof of Corollary 3.8 is made with the help of Theorem 3.7 and system (3.25). Exponential stability follows directly from the fact that for every $(t, x_0, d) \in \Re^+ \times \Re^n \times M_D$ the unique solution $x(t, x_0; d)$ of (1.1) is related to the unique solution $y(t)$ of (3.25) with initial condition $y(0) = x_0$ corresponding to the same $d \in M_D$ by the equation $x(t, x_0; d) = y\left(\int_0^t \phi(x(\tau, x_0; d)) d\tau\right)$. ◁

Since the estimation of the function $\gamma \in K_\infty$ is crucial for the verification of inequalities (3.6), (3.8) and (3.9) of Theorem 3.1 and Corollary 3.5, less conservative estimates of the solution of system (1.1) can be useful. The following theorem uses an additional differential inequality, which can be used to replace inequality (3.8) by a less demanding inequality.

**Corollary 3.9:** *Consider system (1.1) under hypotheses (H1), (H2), (H4) and suppose that there exist a positive definite and radially unbounded function $V \in C^1(\Re^n; \Re^+)$, a locally Lipschitz function $\phi: \Re^n \to (0, +\infty)$, a family of functions $W_i \in C^1(\Re^n; \Re)$ with $W_i(0) = 0$, $b_i \in C^0(\Re^+; \Re^+)$ $(i = 0, ..., k)$, $\rho, c_1, c_2, g \in K_\infty$, with $\rho(s) > c_1(s) \geq c_2(s)$ for all $s > 0$, $\lambda \in K$ with $\lambda(s) < s$ for all $s > 0$, a locally bounded function $r: \Re^+ \to (0, +\infty)$ and a $C^1$ function $\mu: \Re^+ \to \Re$ with $\mu(0) = 0$ for which the function $\kappa(s) := c_1(s) + \mu(s)$ is non-decreasing, such that inequalities (3.22), (3.23), (3.24), (3.3), (3.7) hold. Moreover, suppose that there exist functions $\widetilde{g}, \gamma \in K_\infty$ such that inequalities (3.6), (3.9) hold as well as the following inequalities:*

$$\max_{d \in D} \nabla W_k(x) f(d, x) \leq -\phi(x) \widetilde{g}(V(x)),$$

for all $x \in \Re^n$ with $W_0(x) \geq c_2(V(x))$ and $\max_{i=1,...,k} W_i(x) \geq 0$ \hfill (3.38)

$$\gamma(s) \geq \min \left\{ s, \rho^{-1} \left( \max_{\substack{\tau \in [0, r(s)], \\ V(x) = s}} \left[ \sum_{i=0}^{k} \frac{\tau^i}{i!} W_i(x) - \widetilde{g}(\lambda(s)) \frac{\tau^{k+1}}{(k+1)!} \right] \right) \right\}, \text{ for all } s > 0 \hfill (3.39)$$

*Then $0 \in \Re^n$ is URGAS for (1.1).*

The reader should notice that in general the function $\widetilde{g} \in K_\infty$ involved in (3.38) will be greater than the function $g \in K_\infty$ involved in (3.24). Therefore, (3.39) is a less demanding inequality than (3.8).

**Proof:** It suffices to show that the result holds for the special case $\phi(x) \equiv 1$. Then a similar argument to the one used in the proof of Corollary 3.5 can show the validity of the result to the general case. Therefore, we assume that inequalities (3.1), (3.2), (3.4) and (3.38) with $\phi(x) \equiv 1$ hold.

The reader should notice that inequality (3.8) in the proofs of Theorem 3.1 and Lemma 3.4 is used only for the derivation of inequality $V(t_1) \leq \gamma(V(t_0))$, where $t_0 < t_1$ are times with

$$t_1 - t_0 \leq r(V(t_0)) \hfill (3.40)$$

$$W_0(x(t, x_0; d)) \geq c_2(V(x(t, x_0; d))) \text{ for all } t \in [t_0, t_1] \hfill (3.41)$$

$$V(t) > \lambda(V(t_0)) \text{ for all } t \in [t_0, t_1] \hfill (3.42)$$

Particularly, for Theorem 3.1 we have $t_0 = 0$. Using inequalities (3.2), (3.4), (3.41), (3.42), we obtain that inequalities (A8), (A9) hold for $t \in [t_0, t_1]$ a.e.. Moreover, inequalities (A8), (A9) show that if there exists $T \in [t_0, t_1]$ such that $\max_{i=1,...,k} W_i(T) \leq 0$ then we have $\max_{i=1,...,k} W_i(t) \leq 0$ and $W_0(t) \leq W_0(T)$ for all $t \in [T, t_1]$.

We next distinguish the following cases:

Case 1: $\max_{i=1,...,k} W_i(t) > 0$ for all $t \in [t_0, t_1]$. In this case, inequality (3.38) implies that inequalities (A9), (A10) hold with $g \in K_\infty$ replaced by $\widetilde{g} \in K_\infty$ for all $t \in [t_0, t_1]$ (inequality (A9) holds for $t \in [t_0, t_1]$ a.e.). Consequently, by virtue of (3.40) we get:

$$\max_{t \in [t_0, t_1]} W_0(x(t, x_0; d)) \leq \max_{\tau \in [0, r(V(t_0))]} \left[ \sum_{i=0}^{k} \frac{\tau^i}{i!} W_i(t_0) - \widetilde{g}(\lambda(V(t_0))) \frac{\tau^{k+1}}{(k+1)!} \right] \hfill (3.43)$$

Case 2: There exists $T \in (t_0, t_1]$ such that $\max_{i=1,...,k} W_i(T) \leq 0$. In this case, we must have $\max_{i=1,...,k} W_i(t) \geq 0$ for all $t \in [t_0, T]$. Therefore, inequality (3.38) implies that inequalities (A9), (A10) hold with $g \in K_\infty$ replaced by $\widetilde{g} \in K_\infty$ for all $t \in [t_0, T]$ (inequality (A9) holds for $t \in [t_0, T]$ a.e.). Moreover, since $W_0(t) \leq W_0(T)$ for all $t \in [T, t_1]$, it follows that $\max_{t \in [t_0, t_1]} W_0(x(t, x_0; d)) = \max_{t \in [t_0, T]} W_0(x(t, x_0; d))$. We conclude that (3.43) holds in this case as well.

Case 3: $\max_{i=1,...,k} W_i(t_0) \leq 0$. In this case, we have $W_0(t) \leq W_0(t_0)$ for all $t \in [t_0, t_1]$ and we conclude that inequality (3.43) holds in this case as well.

By virtue of (3.1) and (3.43) we obtain for $t \in [t_0, t_1]$ a.e.:

$$\dot{V}(t) \leq -\rho(V(t)) + \max_{\tau \in [0, r(V(t_0))]} \left[ \sum_{i=0}^{k} \frac{\tau^i}{i!} W_i(t_0) - \tilde{g}(\lambda(V(t_0))) \frac{\tau^{k+1}}{(k+1)!} \right] \tag{3.44}$$

Differential inequality (3.44) directly implies that:

$$V(t) \leq \max \left\{ V(t_0), \rho^{-1} \left( \max_{\tau \in [0, r(V(t_0))]} \left[ \sum_{i=0}^{k} \frac{\tau^i}{i!} W_i(t_0) - \tilde{g}(\lambda(V(t_0))) \frac{\tau^{k+1}}{(k+1)!} \right] \right) \right\}, \text{ for all } t \in [t_0, t_1]$$

The above inequality in conjunction with inequality (3.39) implies that inequality $V(t_1) \leq \gamma(V(t_0))$ holds. It should be noticed that inequality $V(t_1) \leq \gamma(V(t_0))$ holds as well for the case $t_1 = t_0$ (since (3.39) implies that $V(t_0) \leq \gamma(V(t_0))$). The proof is complete. ◁

Similarly with Corollary 3.9, we obtain the following result.

**Corollary 3.10:** *Consider system (1.1) under hypotheses (H1), (H2), (H3) and suppose that there exist a positive definite and radially unbounded function $V \in C^1(\Re^n; \Re^+)$, a locally Lipschitz function $\phi : \Re^n \to (0, +\infty)$, a family of functions $W_i \in C^1(\Re^n; \Re)$ with $W_i(0) = 0$, constants $b_i \geq 0$ ($i = 0,...,k$) with $b_0 \geq \rho$, $\rho, c_1, c_2, g, r > 0$, with $\rho > c_1 \geq c_2 > 0$, $\lambda \in (0,1)$ and $\mu \geq -c_1$ such that inequalities (3.35), (3.36), (3.37), (3.28), (3.30) and (3.32) hold. Moreover, suppose that there exist constants $\tilde{g}, \gamma > 0$ such that inequality (3.31) holds as well as the following inequalities:*

$$\max_{d \in D} \nabla W_k(x) f(d, x) \leq -\tilde{g} \phi(x) V(x),$$

for all $x \in \Re^n$ with $W_0(x) \geq c_2 V(x)$ and $\max_{i=1,...,k} W_i(x) \geq 0$ \qquad (3.45)

$$\gamma \geq \min \left\{ \exp((b_0 - \rho)r), \sup_{\tau \in [0,r], x \in \Re^n} \frac{1}{\rho} \left[ \sum_{i=0}^{k} \frac{\tau^i}{i!} \frac{W_i(x)}{V(x)} - \tilde{g} \lambda \frac{\tau^{k+1}}{(k+1)!} \right] \right\} \tag{3.46}$$

*Then $0 \in \Re^n$ is URGAS for (1.1). Moreover, if there exist constants $0 < K_1 < K_2$ with $K_1 |x|^2 \leq V(x) \leq K_2 |x|^2$ and $\phi(x) \geq K_1$ for all $x \in \Re^n$ then $0 \in \Re^n$ is Robustly Globally Exponentially Stable for (1.1).*

## 4. Examples

This section is devoted to the presentation of two illustrative examples. Both examples can be handled easily by classical Lyapunov analysis (i.e., it is easy to find a continuously differentiable, positive definite and radially unbounded function with negative definite derivative). However, here the issue is to show how we can prove robust global asymptotic (or exponential) stability by using a positive definite function with non sign definite derivative. In both examples, the simplest continuously differentiable, positive definite function $V(x) = |x|^2$ is used; this function fails to satisfy the requirements of Lyapunov's direct method.

**Example 4.1:** Consider the planar system:

$$\dot{x}_1 = -x_1$$
$$\dot{x}_2 = d\beta(x_1) - x_2 \qquad (4.1)$$
$$x = (x_1, x_2)' \in \Re^2, d \in [-p, p]$$

where $p \geq 0$ is a constant parameter, $\beta : \Re \to \Re$ is a locally Lipschitz mapping with $\beta(0) = 0$ and the Lyapunov function

$$V(x) = x_1^2 + x_2^2 \qquad (4.2)$$

The reader should notice that the derivative of the Lyapunov function defined by (4.2) is not necessarily sign definite and classical Lyapunov analysis does not help. Of course there are Lyapunov functions that can be used directly for classical Lyapunov analysis (e.g., $V(x) = x_1^2 + x_2^2 + p^2 \int_0^{x_1} \frac{\beta^2(y)}{y} dy$). Here, for illustration purposes, we apply the result of Theorem 3.1 and we show that for every $p \geq 0$, $0 \in \Re^2$ is URGAS for system (4.1).

We have:

$$\dot{V}(x) = -2V(x) + 2d\beta(x_1)x_2 \leq -V(x) + p^2 \beta^2(x_1), \; \forall x \in \Re^2 \qquad (4.3)$$

Let $\tilde{\beta} : \Re \to \Re$ be an odd $C^1$ mapping such that its restriction on $\Re^+$ is a convex $K_\infty$ function and satisfies $|\beta(x_1)| \leq \tilde{\beta}(|x_1|)$ for all $x_1 \in \Re$. Inequality (4.3) shows that inequality (3.1) holds with $\rho(s) := s$ and $W_0(x) := p^2 \tilde{\beta}^2(x_1)$.

Let $0 < c_1 < 1$ and define $c_1(s) := c_1 s$. Moreover, notice that notice that since $\tilde{\beta} : \Re \to \Re$ is a locally Lipschitz function with $\tilde{\beta}(0) = 0$, there exists $b_0 \in K_\infty \cap C^1(\Re^+; \Re^+)$ such that inequality (3.3) for $i = 0$ holds for all $x \in \Re^2$. Without loss of generality we may assume that $b_0(s) \geq s$ for all $s \geq 0$. Furthermore, since $\tilde{\beta} : \Re \to \Re$ is a $C^1$ mapping, we have:

$$\dot{W}_0(x) = -2p^2 \tilde{\beta}(x_1) \tilde{\beta}'(x_1) x_1, \; \forall x \in \Re^2 \qquad (4.4)$$

If $x_1 \geq 0$, then since the restriction of $\tilde{\beta} : \Re \to \Re$ on $\Re^+$ is a convex $K_\infty$ function we get $\tilde{\beta}(x_1) \leq \tilde{\beta}'(x_1) x_1$. Therefore, (4.4) implies for all $x_1 \geq 0$:

$$\dot{W}_0(x) \leq -2W_0(x) \qquad (4.5)$$

Since $\tilde{\beta} : \Re \to \Re$ is an odd mapping, we have $\tilde{\beta}'(x_1) = \tilde{\beta}(-x_1)$ for all $x_1 < 0$. Therefore, if $x_1 < 0$, by virtue of convexity we get $-\tilde{\beta}(x_1) \leq -\tilde{\beta}'(x_1) x_1$. Thus (4.5) holds for $x_1 < 0$ as well.

Inequality (4.5) shows that inequalities (3.4), (3.5) hold with $k = 0$, $\mu(s) \equiv 0$ and $g(s) := 2c_2(s)$, where $c_2 \in K_\infty$ is an arbitrary function (yet to be selected) that satisfies $c_2(s) \leq c_1 s$ for all $s \geq 0$. Consequently, since $k = 0$, inequalities (3.2) do not apply in this case. Furthermore, since $k = 0$, it follows that inequality (3.8) holds with $\gamma(s) := b_0(s)$ for every locally bounded mapping $r : \Re^+ \to (0, +\infty)$.

Define $\lambda(s) := \lambda s$, where $\lambda \in (0,1)$. Inequality (3.9) holds since $\limsup_{s \to 0^+} \int_{\lambda(s)}^{\gamma(s)} \frac{d\tau}{\rho(\tau) - c_1(\tau)} \leq \frac{1}{1 - c_1} \ln\left(\frac{R}{\lambda}\right) < +\infty$, where $R > 0$ is an appropriate constant that satisfies $b_0(s) \leq R s$ for all $s > 0$ sufficiently small (there exists such a constant since $b_0 \in K_\infty \cap C^1(\Re^+; \Re^+)$).

Finally, notice that inequality (3.6) holds with $c_2(s) := c_1 \lambda^2 b_0^{-1}(s)$. Notice that since $b_0(s) \geq s$ for all $s \geq 0$, we indeed obtain that $c_2(s) \leq \lambda^2 c_1(s) < c_1(s)$ for all $s > 0$. Moreover, since $b_0(s) \leq R s$ for all $s > 0$ sufficiently small, we obtain $s \leq \frac{R}{c_1 \lambda^2} c_2(s)$ for all $s > 0$ sufficiently small.

Thus we are left with inequality (3.7). By virtue of Remark 3.2 and since $b_0(s) \leq R s \leq \frac{R^2}{c_1 \lambda^3} c_2(\lambda s)$ for all $s > 0$ sufficiently small, we get $\limsup_{s \to 0^+} \frac{b_0(s) - c_2(\lambda(s))}{g(\lambda(s))} < +\infty$. Consequently, the locally bounded mapping defined by $r(s) := \frac{1}{2} + \frac{b_0(s)}{2 c_2(\lambda(s))}$ for $s > 0$ and $r(0) := 1$ satisfies inequality (3.7).

Therefore, all requirements of Theorem 3.1 hold. We conclude that for every $p \geq 0$, $0 \in \Re^2$ is URGAS for system (4.1).

If we further assume that $|\beta(x_1)| \leq K|x_1|$, where $K > 0$ then the reader can verify that all requirements of Theorem 3.7 hold. Particularly, all the above hold with $b_0(s) := (1 + K^2)s$. Indeed, in this case we may conclude that for every $p \geq 0$, $0 \in \Re^2$ is URGES for system (4.1). ◁

**Example 4.2:** Consider the linear uncertain system

$$\begin{aligned} \dot{x}_1 &= x_2 \\ \dot{x}_2 &= -(1+d)x_1 - 2x_2 \\ x &= (x_1, x_2)' \in \Re^2, d \in D := [0, p] \end{aligned} \quad (4.6)$$

where $p \geq 0$ is a constant parameter. Our goal is to determine the maximum allowable value of $p \geq 0$ for which $0 \in \Re^2$ is Robustly Globally Exponentially Stable. To this purpose, we will use the Lyapunov function defined by (4.2) and Corollary 3.10 with $\phi(x) \equiv 1$. It should be noticed that the derivative of $V$ is not negative definite (it is only negative semi-definite only for the case $d \equiv 0$). Indeed, we have by completing the squares

$$\max_{d \in D} \nabla V(x) f(d, x) = \max_{d \in D} \left( -2 d x_1 x_2 - 4 x_2^2 \right) \leq p^2 x_1^2 - 3 x_2^2 = -3V(x) + (3 + p^2) x_1^2 \quad (4.7)$$

where $f(d, x) := (x_2, -(1+d)x_1 - 2x_2)'$. Inequality (4.7) shows that inequality (3.26) holds with $\rho := 3$ and $W_0(x) := (3 + p^2) x_1^2$. We also have:

$$\max_{d \in D} \nabla W_0(x) f(d, x) = 2(3 + p^2) x_1 x_2 \quad (4.8)$$

Equation (4.8) shows that inequality (3.27) for $i = 0$ holds for arbitrary $c_2 > 0$ with $W_1(x) := 2(3 + p^2) x_1 x_2$. It should be noticed that the inequality $W_0(x) \geq c_2 V(x)$ is equivalent to the inequality

$$-\sqrt{\frac{3 + p^2 - c_2}{c_2}} \leq \frac{x_2}{x_1} \leq \sqrt{\frac{3 + p^2 - c_2}{c_2}} \quad (4.9)$$

for every $c_2 > 0$. In addition, we get

$$\max_{d \in D} \nabla W_1(x) f(d, x) = 2(3 + p^2) \max_{d \in D} \left( x_2^2 - (1+d) x_1^2 - 2 x_1 x_2 \right) \leq 2(3 + p^2)\left( x_2^2 - x_1^2 - 2 x_1 x_2 \right) \quad (4.10)$$

Inequality (4.10) shows that there exists a constant $g > 0$ such that $\max_{d \in D} \nabla W_1(x) f(d,x) \leq -2(3+p^2)gV(x)$, provided that

$$\frac{1-\sqrt{2-g^2}}{1+g} \leq \frac{x_2}{x_1} \leq \frac{1+\sqrt{2-g^2}}{1+g} \tag{4.11}$$

By virtue of (4.9) and (4.10), it follows that inequality (3.29) for $k=1$ holds as long as

$$\frac{3+p^2}{4-2\sqrt{2}} < c_2 \tag{4.12}$$

Since the inequality $c_2 < \rho = 3$ must also hold, it follows from (4.12) that the maximum allowable value of $p \geq 0$ must satisfy $p < \sqrt{9 - 6\sqrt{2}}$. In order to determine the constant $\tilde{g} > 0$ that satisfies (3.45) with $\phi(x) \equiv 1$, we notice that inequalities (4.9) and (4.10) in conjunction with the constraint $W_1(x) \geq 0$ and (4.12) (which implies $\sqrt{\frac{3+p^2-c_2}{c_2}} < \sqrt{2}-1$) give:

$$\tilde{g} = 2(3+p^2) \min\left\{ \frac{1+2y-y^2}{1+y^2} : 0 \leq y \leq \sqrt{\frac{3+p^2-c_2}{c_2}} \right\} = 2(3+p^2) \tag{4.13}$$

The reader should notice that inequalities (3.28) hold with $b_0 = b_1 = 3+p^2$. Moreover, inequality (3.32) holds with

$$r = 1 + 2\sqrt{\frac{3+p^2}{\lambda g}} \max\left\{1, 2\sqrt{\frac{3+p^2}{\lambda g}}\right\} \text{ (see Remark 3.2) for every } \lambda \in (0,1).$$

We next determine constants $\mu \geq 0$ and $c_1 \in [c_2, 3)$ so that inequality (3.230) holds. We have

$$\max_{d \in D} \nabla W_0(x) f(d,x) + \mu \max_{d \in D} \nabla V(x) f(d,x) \leq 2(3+p^2)x_1 x_2 + \mu p^2 x_1^2 - 3\mu x_2^2 \tag{4.14}$$

The reader can verify that $c_2 V(x) \leq W_0(x) \leq c_1 V(x)$ is equivalent to the sector condition:

$$\sqrt{\frac{3+p^2-c_2}{c_2}} \geq \frac{x_2}{x_1} \geq \sqrt{\frac{3+p^2-c_1}{c_1}} \quad \text{or} \quad -\sqrt{\frac{3+p^2-c_2}{c_2}} \leq \frac{x_2}{x_1} \leq -\sqrt{\frac{3+p^2-c_1}{c_1}} \tag{4.15}$$

On the other hand, the right hand side of inequality (4.14) is non-positive provided that

$$\frac{x_2}{x_1} \geq \frac{1+\sqrt{1+12s^2 p^2}}{6s} \quad \text{or} \quad \frac{x_2}{x_1} \leq \frac{1-\sqrt{1+12s^2 p^2}}{6s} \tag{4.16}$$

where $s = \mu / 2(3+p^2)$. It follows that (4.15) implies (4.16) if the following inequality holds:

$$\mu \geq \frac{2\sqrt{c_1(3+p^2-c_1)}}{3-c_1} \tag{4.17}$$

Thus we are left with the verification of inequalities (3.31) and (3.46). By virtue of (4.13) and previous definitions the following inequalities hold for every $\lambda \in (0,1)$:

$$\sup_{\tau\in[0,r],x\in\Re^n}\frac{1}{\rho}\left[\sum_{i=0}^{k}\frac{\tau^i}{i!}\frac{W_i(x)}{V(x)}-\widetilde{g}\lambda\frac{\tau^{k+1}}{(k+1)!}\right]\leq\sup_{\tau\in[0,r],x\in\Re^n}\frac{3+p^2}{3}\left[\frac{x_1^2}{V(x)}+\tau\frac{2x_1x_2}{V(x)}-\lambda\tau^2\right]\leq$$

$$\leq\sup_{x\in\Re^n}\frac{3+p^2}{3}\left[\frac{x_1^2}{V(x)}+\frac{x_1^2x_2^2}{\lambda V^2(x)}\right]=\sup_{x\in\Re^n}\frac{3+p^2}{3\lambda}\left[(1+\lambda)\frac{x_1^2}{V(x)}-\frac{x_1^4}{V^2(x)}\right]\leq\frac{3+p^2}{12\lambda}(\lambda+1)^2$$

Consequently, inequality (3.46) will hold with $\gamma:=\frac{3+p^2}{12\lambda}(\lambda+1)^2$. On the other hand, previous definitions imply that inequality (3.31) is equivalent to the following inequality:

$$\frac{12\lambda^2}{(\lambda+1)^2}\frac{c_1(3-c_1)+2\sqrt{c_1(3+p^2-c_1)}}{c_2(3-c_1)+2\sqrt{c_1(3+p^2-c_1)}}-3>p^2$$

for arbitrary constants $\frac{3+p^2}{4-2\sqrt{2}}<c_2<c_1<3$ and $\lambda\in(0,1)$. Consequently, the maximum allowable value of $p\geq 0$ must satisfy $p<\sqrt{9-6\sqrt{2}}$ and the following inequality

$$\frac{3(c_1-c_2)(3-c_1)}{c_2(3-c_1)+2\sqrt{c_1(3+p^2-c_1)}}>p^2 \tag{4.18}$$

for certain constants $\frac{3+p^2}{4-2\sqrt{2}}<c_2<c_1<3$. Numerical calculations show that the maximum value is greater than $\frac{1}{5}\sqrt{\frac{7}{5}}\approx 0.236643$; the reader can verify that inequality (4.18) holds with $p=\frac{1}{5}\sqrt{\frac{7}{5}}$, $c_2=2.6094$ and $c_1=2.8594$.

It should be noticed that the result is very conservative. Indeed, by following classical Lyapunov analysis the reader can verify that much higher values for $p\geq 0$ than $\frac{1}{5}\sqrt{\frac{7}{5}}\approx 0.236643$ can be allowed. For example, the quadratic Lyapunov function $V(x)=\frac{1}{4}x_1^2+\frac{1}{2}(x_2+\sigma x_1)^2$ with $\sigma=1+\frac{\sqrt{2}}{2}$ has negative definite derivative for $p<1$. However, here we have used a completely inappropriate Lyapunov function, which has positive derivative in certain regions of the state space. The example simply shows that stability analysis is possible even with completely inappropriate Lyapunov functions. ◁

## 5. Concluding Remarks

Novel criteria for global asymptotic stability are presented. The results (Theorem 3.1, Corollary 3.5, Theorem 3.7, Corollaries 3.8, 3.9 and 3.10) are developed for the autonomous uncertain case and are obtained by a combination of:

- suitable generalizations of the "discretization approach" (Proposition 2.3 and Proposition 2.5), which are necessary and sufficient conditions for uniform robust global asymptotic stability,

- the idea contained in the proof of the original Matrosov's result concerning the division of the state space into two regions: the "good region", where the derivative of the Lyapunov function has a negative upper bound and the "bad region" where the derivative of the Lyapunov function can be positive.

The results can be used for the proof of global asymptotic stability by using continuously differentiable, positive definite functions which do not have a negative semi-definite derivative. Illustrating examples are provided, which show how we can use very simple positive definite functions (e.g. $V(x)=|x|^2$), which do not have a sign-definite derivative.

Future work can address the issue of the extension of the obtained results to the local case or the time-varying case.

**Acknowledgments:** The author would like to thank Professor A. R. Teel and Professor J. Tsinias for their useful comments and suggestions.## References

[1] Aeyels, D. and J. Peuteman, "A New Asymptotic Stability Criterion for Nonlinear Time-Variant Differential Equations", *IEEE Transactions on Automatic Control*, 43(7), 1998, 968-971.

[2] Coron, J.-M. and L. Rosier, "A Relation Between Continuous Time-Varying and Discontinuous Feedback Stabilization", *Journal of Mathematical Systems, Estimation, and Control*, 4(1), 1994, 67-84.

[3] Fillipov, A. V., "Differential Equations with Discontinuous Right-Hand Sides", Kluwer Academic Publishers, 1988.

[4] Jiang, Z.-P. and Y. Wang, "A Converse Lyapunov Theorem for Discrete-Time Systems with Disturbances", *Systems and Control Letters*, 45, 2002, 49-58.

[5] Karafyllis, I., "Non-Uniform in Time Robust Global Asymptotic Output Stability", *Systems and Control Letters*, 54(3), 2005, 181-193.

[6] Karafyllis, I., "A System-Theoretic Framework for a Wide Class of Systems II: Input-to-Output Stability", *Journal of Mathematical Analysis and Applications*, 328(1), 2007, 466-486.

[7] Karafyllis, I., and J. Tsinias, "Control Lyapunov Functions and Stabilization by Means of Continuous Time-Varying Feedback", *ESAIM Control, Optimisation and Calculus of Variations*, 15(3), 2009, 599-625.

[8] Karafyllis, I., C. Kravaris and N. Kalogerakis, "Relaxed Lyapunov Criteria for Robust Global Stabilization of Nonlinear Systems", to appear in the *International Journal of Control*.

[9] Khalil, H.K., *Nonlinear Systems*, 2$^{nd}$ Edition, Prentice-Hall, 1996.

[10] Lin, Y., E.D. Sontag and Y. Wang, "A Smooth Converse Lyapunov Theorem for Robust Stability", *SIAM Journal on Control and Optimization*, 34, 1996, 124-160.

[11] Loria, A., E. Panteley, D. Popovic and A. R. Teel, "A Nested Matrosov Theorem and Persistency of Excitation for Uniform Convergence in Stable Nonautonomous Systems", *IEEE Transactions on Automatic Control*, 50(2), 2005, 183-198.

[12] Malisoff, M. and F. Mazenc, "Constructions of strict Lyapunov functions for discrete time and hybrid time-varying systems", *Nonlinear Analysis: Hybrid Systems*, 2(2), 2008, 394-407.

[13] Mazenc, F. and D. Nesic, "Strong Lyapunov Functions for Systems Satisfying the Conditions of La Salle", *IEEE Transactions on Automatic Control*, 49(6), 2004, 1026-1030.

[14] Mazenc, F. and D. Nesic, "Lyapunov Functions for Time-Varying Systems Satisfying Generalized Conditions of Matrosov Theorem", *Mathematics of Control, Signals and Systems*, 19, 2007, 151-182.

[15] Mazenc, F., M. Malisoff, and O. Bernard, "A simplified design for strict Lyapunov functions under Matrosov conditions", *IEEE Transactions on Automatic Control*, 54(1), 2009, 177-183.

[16] Malisoff, M., and F. Mazenc, *Constructions of Strict Lyapunov Functions*, Communications and Control Engineering Series, Springer-Verlag, London, 2009.

[17] Michel, A. N., L. Hou and D. Liu, *Stability of Dynamical Systems Continuous, Discontinuous and Discrete Systems*, Birkhauser, Boston, 2008.

[18] Michel, A. N., and L. Hou, "Stability Results Involving Time-Averaged Lyapunov Function Derivatives", *Nonlinear Analysis: Hybrid Systems*, 3, 2009, 51-64.

[19] Munoz de la Pena, D., and P. D. Christofides, "Stability of Nonlinear Asynchronous Systems", *Systems and Control Letters*, 57, 2008, 465-473.

[20] Peuteman, J. and D. Aeyels, "Exponential Stability of Slowly Time-Varying Nonlinear Systems", *Mathematics of Control, Signals and Systems*, 15, 2002, 42-70.

[21] Peuteman, J. and D. Aeyels, "Exponential Stability of Nonlinear Time-Varying Differential Equations and Partial Averaging", *Mathematics of Control, Signals and Systems*, 15, 2002, 202-228.

[22] Rouche, N., P. Habets and M. Laloy, *Stability Theory by Liapunov's Direct Method*, Springer-Verlag, New York, 1977.

[23] Stuart, A. M. and A. R. Humphries, *Dynamical Systems and Numerical Analysis*, Cambridge University Press, 1998.

[24] Teel, A. R., E. Panteley and A. Loria, "Integral Characterizations of Uniform Asymptotic and Exponential Stability with Applications", *Mathematics of Control, Signals and Systems*, 15, 2002, 177-201.19

# Appendix

**Proof of FACT I in the proof of Theorem 3.1:** Suppose on the contrary that $W_0(x(t, x_0; d)) \geq c_2(V(x(t, x_0; d)))$ for all $t \in [0, r(V(x_0))]$. Using inequalities (3.2), (3.4) and the fact that $V(x(t, x_0; d)) > \lambda(V(x_0))$ for all $t \in [0, r(V(x_0))]$, we obtain for $t \in [0, r(V(x_0))]$ a.e.:

$$\dot{W}_i(t) \leq W_{i+1}(t), \text{ for } i = 0, ..., k-1 \tag{A1}$$

$$\dot{W}_k(t) \leq -g(\lambda(V(x_0))) \tag{A2}$$

where $W_i(t) := W_i(x(t, x_0; d))$ ($i = 0, ..., k$). Inequalities (A1) and (A2) imply that the following inequality holds for all $t \in [0, r(V(x_0))]$

$$W_0(t) \leq \sum_{i=0}^{k} \frac{t^i}{i!} W_i(x_0) - g(\lambda(V(x_0))) \frac{t^{k+1}}{(k+1)!} \tag{A3}$$

Our assumption that $W_0(x(t, x_0; d)) \geq c_2(V(x(t, x_0; d)))$ for all $t \in [0, r(V(x_0))]$ in conjunction with the fact that $V(x(t, x_0; d)) > \lambda(V(x_0))$ for all $t \in [0, r(V(x_0))]$ and inequality (A3) gives:

$$c_2(\lambda(V(x_0))) \leq \sum_{i=0}^{k} \frac{t^i}{i!} W_i(x_0) - g(\lambda(V(x_0))) \frac{t^{k+1}}{(k+1)!}, \text{ for all } t \in [0, r(V(x_0))] \tag{A4}$$

Inequality (A4) (for $t = r(V(x_0))$) in conjunction with inequalities (3.3) (which give $W_i(x_0) \leq b_i(V(x_0))$ for $i = 0, ..., k$) implies that the following inequality must hold:

$$c_2(\lambda(V(x_0))) \leq \sum_{i=0}^{k} \frac{r^i(V(x_0))}{i!} b_i(V(x_0)) - g(\lambda(V(x_0))) \frac{r^{k+1}(V(x_0))}{(k+1)!}$$

which contradicts inequality (3.7) with $s = V(x_0)$. The proof is complete. ◁

**Proof of FACT II in the proof of Theorem 3.1:** Suppose first that $t_0 > 0$. By virtue of definition (3.13), it follows that $W_0(x(t, x_0; d)) \geq c_2(V(x(t, x_0; d)))$ for all $t \in [0, t_0]$. Using inequalities (3.2), (3.4) and the fact that $V(x(t, x_0; d)) > \lambda(V(x_0))$ for all $t \in [0, r(V(x_0))]$, it follows that inequalities (A1), (A2) hold for $t \in [0, t_0]$ a.e.. Consequently, inequality (A3) holds for all $t \in [0, t_0]$. Using inequalities (3.3) (which give $W_i(x_0) \leq b_i(V(x_0))$ for $i = 0, ..., k$), we get from (A3):

$$W_0(t) \leq \sum_{i=0}^{k} \frac{t^i}{i!} b_i(V(x_0)) - g(\lambda(V(x_0))) \frac{t^{k+1}}{(k+1)!} \tag{A5}$$

Inequality (A5) and the fact that $t_0 \in [0, r(V(x_0))]$ imply that:

$$\max_{t \in [0, t_0]} W_0(x(t, x_0; d)) \leq \max_{t \in [0, r(V(x_0))]} \left[ \sum_{i=0}^{k} \frac{t^i}{i!} b_i(V(x_0)) - g(\lambda(V(x_0))) \frac{t^{k+1}}{(k+1)!} \right] \tag{A6}$$

By virtue of (3.1) and (A6) we obtain for $t \in [0, t_0]$ a.e.:

$$\dot{V}(t) \leq -\rho(V(t)) + \max_{\tau \in [0, r(V(x_0))]} \left[ \sum_{i=0}^{k} \frac{\tau^i}{i!} b_i(V(x_0)) - g(\lambda(V(x_0))) \frac{\tau^{k+1}}{(k+1)!} \right] \tag{A7}$$

where $V(t) := V(x(t, x_0; d))$. Differential inequality (A7) directly implies that:

$$V(t) \leq \max\left\{ V(x_0), \rho^{-1}\left( \max_{\tau \in [0, r(V(x_0))]} \left[ \sum_{i=0}^{k} \frac{\tau^i}{i!} b_i(V(x_0)) - g(\lambda(V(x_0))) \frac{\tau^{k+1}}{(k+1)!} \right] \right) \right\}, \text{ for all } t \in [0, t_0]$$

The above inequality in conjunction with inequality (3.8) implies that inequality (3.14) holds. It should be noticed that inequality (3.14) holds as well for the case $t_0 = 0$ (since (3.8) implies that $V(x_0) \leq \gamma(V(x_0))$).

The proof is complete.  ◁

**Proof of Lemma 3.4:** We will prove that system (1.1) is forward complete by contradiction. Suppose that there exists $x_0 \in \Re^n \setminus \{0\}$, $d \in M_D$ the solution $x(t, x_0; d)$ of (1.1) is defined on $[0, t_{\max})$, where $t_{\max} \in (0, +\infty)$ and cannot be further continued. Standard results on the continuation of the solutions of ordinary differential equations imply that $\lim_{t \to t_{\max}^-} V(t) = +\infty$, where $V(t) := V(x(t, x_0; d))$.

We next prove the following claims.

CLAIM 1: There exists $t_0 \in [0, t_{\max})$ such that $V(t) > \lambda(V(t_0))$, for all $t \in [t_0, t_{\max})$.

Proof of Claim 1: If Claim 1 were not true then for every $t_i \in [0, t_{\max})$, there would exist $t_{i+1} \in (t_i, t_{\max})$ with $V(t_{i+1}) \leq \lambda(V(t_i))$. Consequently, we can construct an increasing sequence $\{t_i\}_{i=0}^{\infty}$ with $V(t_{i+1}) \leq \lambda(V(t_i))$, $V(t_i) \leq \lambda^{(i)}(V(t_0))$, where $\lambda^{(i)}(s) := \underbrace{(\lambda \circ \lambda \circ \ldots \circ \lambda)}_{i \text{ times}}(s)$ for $i \geq 1$ and $\lambda^{(0)}(s) := s$. A standard contradiction argument shows that $t_i \to T := \sup t_i \leq t_{\max}$ and $V(t_i) \to 0$.

If $T = t_{\max}$ then we obtain $\liminf_{t \to t_{\max}^-} V(t) = 0$, a contradiction with the fact that $\lim_{t \to t_{\max}^-} V(t) = +\infty$.

If $T < t_{\max}$ then we must have $\liminf_{t \to T^-} V(t) = 0$. Since the mapping $t \to V(t)$ is continuous, we must have $\liminf_{t \to T^-} V(t) = \lim_{t \to T} V(t) = V(T)$ and this implies $V(T) = 0$. Consequently, we must have $x(T, x_0; d) = 0$. Uniqueness of solutions for system ($\Sigma$) implies that $x(t, x_0; d) = 0$, for all $t \geq T$, which contradicts the fact that $\lim_{t \to t_{\max}^-} V(t) = +\infty$.

CLAIM 2: There exists $t \in [t_0, t_{\max})$ with $W_0(x(t, x_0; d)) < c_2(V(x(t, x_0; d)))$.

Proof of Claim 2: Suppose on the contrary that $W_0(x(t, x_0; d)) \geq c_2(V(x(t, x_0; d)))$ for all $t \in [t_0, t_{\max})$. Using inequalities (3.2), (3.4) and the fact that $V(x(t, x_0; d)) > \lambda(V(t_0))$ for all $t \in [t_0, t_{\max})$, we obtain for $t \in [t_0, t_{\max})$ a.e.:

$$\dot{W}_i(t) \leq W_{i+1}(t), \text{ for } i = 0, \ldots, k-1 \tag{A8}$$

$$\dot{W}_k(t) \leq -g(\lambda(V(t_0))) \tag{A9}$$

where $W_i(t) := W_i(x(t, x_0; d))$ ($i = 0, \ldots, k$). Inequalities (A8) and (A9) imply that the following inequality holds for all $t \in [t_0, t_{\max})$

$$W_0(t) \leq \sum_{i=0}^{k} \frac{(t - t_0)^i}{i!} W_i(t_0) - g(\lambda(V(t_0))) \frac{(t - t_0)^{k+1}}{(k+1)!} \tag{A10}$$

Our assumption that $W_0(x(t,x_0;d)) \geq c_2(V(x(t,x_0;d)))$ for all $t \in [t_0, t_{\max})$ in conjunction with inequality (A10) gives:

$$V(t) \leq c_2^{-1}\left(\sum_{i=0}^{k} \frac{(t-t_0)^i}{i!} W_i(t_0) - g(\lambda(V(t_0))) \frac{(t-t_0)^{k+1}}{(k+1)!}\right), \text{ for all } t \in [t_0, t_{\max})$$

which contradicts the fact that $\lim_{t \to t_{\max}^-} V(t) = +\infty$.

Define:

$$t_1 := \inf\{t \in [t_0, t_{\max}) : W_0(x(t, x_0; d)) < c_2(V(x(t, x_0; d)))\} \tag{A11}$$

CLAIM 3: $t_1 - t_0 \leq r(V(t_0))$

Proof of Claim 3: Suppose on the contrary that $t_1 - t_0 > r(V(t_0))$. Then definition (A11) implies that $W_0(x(t, x_0; d)) \geq c_2(V(x(t, x_0; d)))$ for all $t \in [t_0, t_1]$. Using inequalities (3.2), (3.4) and the fact that $V(x(t, x_0; d)) > \lambda(V(t_0))$ for all $t \in [t_0, t_{\max})$, we obtain that inequalities (A8), (A9) hold for $t \in [t_0, t_1]$ a.e.. Inequalities (A8) and (A9) imply that inequality (A10) holds for all $t \in [t_0, t_1]$. Since $W_0(x(t, x_0; d)) \geq c_2(V(x(t, x_0; d)))$ for all $t \in [t_0, t_1]$ and $V(x(t, x_0; d)) > \lambda(V(t_0))$ for all $t \in [t_0, t_{\max})$, we get from (A10):

$$c_2(\lambda(V(t_0))) \leq \sum_{i=0}^{k} \frac{(t-t_0)^i}{i!} W_i(t_0) - g(\lambda(V(t_0))) \frac{(t-t_0)^{k+1}}{(k+1)!}, \text{ for all } t \in [t_0, t_1]$$

Using inequalities (3.3) (which give $W_i(t_0) \leq b_i(V(t_0))$ for $i = 0, \ldots, k$), we get from the above inequality:

$$c_2(\lambda(V(t_0))) \leq \sum_{i=0}^{k} \frac{(t-t_0)^i}{i!} b_i(V(t_0)) - g(\lambda(V(t_0))) \frac{(t-t_0)^{k+1}}{(k+1)!}, \text{ for all } t \in [t_0, t_1]$$

The above inequality for $s = V(t_0)$ and $t = t_0 + r(V(t_0))$ contradicts inequality (3.7).

CLAIM 4: The following inequality holds:

$$V(t_1) \leq \gamma(V(t_0)) \tag{A12}$$

Proof of Claim 4: Suppose first that $t_1 > t_0$. By virtue of definition (A11), it follows that $W_0(x(t, x_0; d)) \geq c_2(V(x(t, x_0; d)))$ for all $t \in [t_0, t_1]$. Using inequalities (3.2), (3.4) and the fact that $V(x(t, x_0; d)) > \lambda(V(t_0))$ for all $t \in [t_0, t_{\max})$, we obtain that inequalities (A8), (A9) hold for $t \in [t_0, t_1]$ a.e.. Inequalities (A8) and (A9) imply that inequality (A10) holds for all $t \in [t_0, t_1]$. Using inequalities (3.3) (which give $W_i(t_0) \leq b_i(V(t_0))$ for $i = 0, \ldots, k$), we get from (A10):

$$W_0(t) \leq \sum_{i=0}^{k} \frac{(t-t_0)^i}{i!} b_i(V(t_0)) - g(\lambda(V(t_0))) \frac{(t-t_0)^{k+1}}{(k+1)!} \tag{A13}$$

Inequality (A13) and Claim 3 imply that:

$$\max_{t \in [t_0, t_1]} W_0(x(t, x_0; d)) \leq \max_{\tau \in [0, r(V(t_0))]} \left[\sum_{i=0}^{k} \frac{\tau^i}{i!} b_i(V(t_0)) - g(\lambda(V(t_0))) \frac{\tau^{k+1}}{(k+1)!}\right] \tag{A14}$$

By virtue of (3.1) and (A14) we obtain for $t \in [t_0, t_1]$ a.e.:

$$\dot{V}(t) \leq -\rho(V(t)) + \max_{\tau \in [0, r(V(t_0))]} \left[ \sum_{i=0}^{k} \frac{\tau^i}{i!} b_i(V(t_0)) - g(\lambda(V(t_0))) \frac{\tau^{k+1}}{(k+1)!} \right] \quad (A15)$$

Differential inequality (A15) directly implies that:

$$V(t) \leq \max \left\{ V(t_0), \rho^{-1}\left( \max_{\tau \in [0, r(V(t_0))]} \left[ \sum_{i=0}^{k} \frac{\tau^i}{i!} b_i(V(t_0)) - g(\lambda(V(t_0))) \frac{\tau^{k+1}}{(k+1)!} \right] \right) \right\}, \text{ for all } t \in [t_0, t_1]$$

The above inequality in conjunction with inequality (3.8) implies that inequality (A12) holds. It should be noticed that inequality (A12) holds as well for the case $t_1 = t_0$ (since (3.8) implies that $V(t_0) \leq \gamma(V(t_0))$).

CLAIM 5: There exists $t \in [t_1, t_{\max})$ with $W_0(x(t, x_0; d)) > c_1(V(x(t, x_0; d)))$.

Proof of Claim 5: Suppose the contrary, that $W_0(x(t, x_0; d)) \leq c_1(V(x(t, x_0; d)))$ for all $t \in [t_1, t_{\max})$. Then inequality (3.1) and the fact that $\rho(s) > c_1(s)$ for all $s > 0$, imply that $\dot{V}(t) \leq 0$ for $t \in [t_1, t_{\max})$ a.e.. Thus we obtain $V(t) \leq V(t_1)$ for all $t \in [t_1, t_{\max})$, a contradiction with the fact that $\lim_{t \to t_{\max}^-} V(t) = +\infty$.

We are now ready to finish the proof. Continuity of mapping $t \to \frac{W_0(x(t, x_0; d))}{V(x(t, x_0; d))}$ guarantees the existence of times $t_2 < t_3$ with $t_1 \leq t_2 < t_3 < t_{\max}$ and such that:

$$W_0(x(t, x_0; d)) \leq c_1(V(x(t, x_0; d))), \text{ for all } t \in [t_1, t_3] \quad (A16)$$

$$W_0(x(t_2, x_0; d)) = c_2(V(x(t_2, x_0; d))), \; W_0(x(t_3, x_0; d)) = c_1(V(x(t_3, x_0; d))) \quad (A17)$$

$$W_0(x(t, x_0; d)) \geq c_2(V(x(t, x_0; d))), \text{ for all } t \in [t_2, t_3] \quad (A18)$$

Inequality (A16) in conjunction with inequality (3.1) guarantees that

$$V(t) \leq V(t_2) - \int_{t_2}^{t} (\rho(V(\tau)) - c_1(V(\tau))) d\tau, \text{ for all } t \in [t_2, t_3] \quad (A19)$$

$$V(t) \leq V(t_1), \text{ for all } t \in [t_1, t_3] \quad (A20)$$

Inequalities (A16), (A18) and (3.5) imply that

$$W_0(t_3) + \mu(V(t_3)) \leq W_0(t_2) + \mu(V(t_2)) \quad (A21)$$

It follows from (A17) and (A21) that:

$$c_1(V(t_3)) + \mu(V(t_3)) \leq c_2(V(t_2)) + \mu(V(t_2)) \quad (A22)$$

If $V(t_3) \leq \lambda(\gamma^{-1}(V(t_2)))$ then using (A12), (A20) we obtain $V(t_3) \leq \lambda(V(t_0))$, a contradiction with Claim 1.

Thus we are left with the case $V(t_3) > \lambda(\gamma^{-1}(V(t_2)))$. In this case inequality (A22) and the fact that the function $\kappa(s) := c_1(s) + \mu(s)$ is non-decreasing, give

$$c_1\left(\lambda\left(\gamma^{-1}(V(t_2))\right)\right) + \mu\left(\lambda\left(\gamma^{-1}(V(t_2))\right)\right) \leq c_2(V(t_2)) + \mu(V(t_2))$$

The above inequality contradicts inequality (3.6) for $s = \gamma^{-1}(V(t_2))$.

The proof is complete. ◁